\documentclass[12pt]{article}
\usepackage{amsmath,amsfonts,amssymb,theorem,graphicx,verbatim,fullpage}
\usepackage{hyperref}
\usepackage{pb-diagram}
\newcommand{\CC}{\mathbb C}

\newcommand{\AAA}{\mathbb A}
\newcommand{\PP}{\mathbb P}
\newcommand{\ZZ}{\mathbb Z}

\newcommand{\QQ}{\mathbb Q}

\newcommand{\lam}{\lambda}

\newcommand{\Lam}{\Lambda}
\newcommand{\al}{\alpha}

\newcommand{\gog}{\mathfrak g}

\newcommand{\gob}{\mathfrak b}
\newcommand{\got}{\mathfrak t}
\newcommand{\gop}{\mathfrak p}
\newcommand{\gogl}{\mathfrak {gl}}
\newcommand{\nabp}{\nabla_{\gop}}

\newcommand{\gosp}{\mathfrak {sp}}
\newcommand{\la}{\lambda}
\newcommand{\om}{\omega}
\newcommand{\si}{\sigma}
\newcommand{\ep}{\varepsilon}
\newcommand{\Si}{\Sigma}

\newcommand{\Oh}{\mathcal O}

\newcommand{\Cc}{\mathcal C}

\newcommand{\into}{\hookrightarrow}
\newcommand{\onto}{\twoheadrightarrow}

\newcommand{\color}[6]{}

\theorembodyfont{\rmfamily}
 \newtheorem{theorem}[subsection]{Theorem}

\newtheorem{example}[subsection]{Example}
\newtheorem{rmk}[subsection]{Remark}
{
\theorembodyfont{\rmfamily}

 \newtheorem{nothing}[subsection]{}
}
 
 \newenvironment{proof}{\paragraph{Proof}}{\par\medskip}

\theorembodyfont{\rmfamily}

\DeclareMathOperator{\fl}{FL}
\DeclareMathOperator{\Sp}{Sp}
\DeclareMathOperator{\codim}{codim}

\DeclareMathOperator{\SL}{SL}
\DeclareMathOperator{\GL}{GL}

\DeclareMathOperator{\Gr}{Gr }
\DeclareMathOperator{\lgr}{LGr }
\DeclareMathOperator{\Hom}{Hom}

\DeclareMathOperator{\Pic}{Pic}
\DeclareMathOperator{\Proj}{Proj}

\DeclareMathOperator{\so}{SO}
\DeclareMathOperator{\ogr}{OGr}

\newtheorem{thm}{Theorem}[section]

 \newtheorem{dfn}[thm]{Definition}

\numberwithin{equation}{section}
\usepackage[normalem]{ulem}
\begin{document}

\title{Constructing projective varieties \\ in weighted flag varieties II }
\date{}
\author{Muhammad Imran Qureshi}
\maketitle

\vspace{0.1in}

\begin{abstract}
We give the construction of weighted Lagrangian Grassmannians \(w\lgr(3,6)\) and weighted partial \(A_3\) flag variety \(w\fl_{1,3}\) coming from the symplectic Lie group \(\Sp(6,\CC)\) and the general linear group \(\GL(4,\CC)\) respectively. We  give  general formulas for their Hilbert series in terms of Lie theoretic  data. We use them as key varieties (Format) to construct some families of polarized 3-folds in codimension 7 and 9. At the end, we list all the distinct weighted flag varieties in codimension \(4\le c\le 10\). \end{abstract}

\section{Introduction}
This article is a sequel of  \cite{qs}.
In that paper, we proved a general formula for the Hilbert series of weighted flag varieties  and discussed a computer aided method to find their defining equations. We also applied them to construct families of projective    varieties  as  quasilinear sections of  certain weighted flag varieties in  codimensions 8 and 6. In \cite{qs}, we used a weighted $G_2$ variety as a key variety(ambient  weighted projective variety)  to construct  families of  polarized varieties in codimension eight and weighted Grassmannian $w\Gr(2,6)$ for the construction of varieties     in  codimension six.

In this paper, we explicitly discuss  the construction   of two new weighted flag varieties \(\left(w\Si,\Oh_{w\Si}(1)\right)\). The first one is the weighted Lagrangian Grassmannian $w\lgr(3,6)$; a homogeneous variety for   the symplectic Lie group $\Sp(6,\CC)$. The second one is the weighted partial flag variety $w\fl_{1,3}$; a partial flag  variety for the general linear group $\GL(4,\CC)$.  The $w\lgr(3,6)$ is a six dimensional variety and has an embedding in weighted projective space $w\PP^{13}$:  a codimension seven embedding. The weighted partial flag variety $w\fl_{1,3}$ is a five dimensional variety and has an embedding in $w\PP^{14}$:  a codimension nine embedding. 

We use these key varieties to exhibit some potentially new  families of Calabi--Yau 3-folds in codimension  seven and nine as weighted complete intersections of these varieties. The constructed Calabi--Yau 3-folds \(X\)  have canonical  
 singularities which can be resolved by    crepant resolutions $Y\to X$ using standard theory \cite{YPG}. The desingularization \(Y\) may lead  to new examples of Calabi--Yau 3-folds but since we do not compute the topological invariants such as   Betti and Hodge numbers, we have no  evidence.

  The explicit constructions of \(w\lgr(3,6)\) and \(w\fl_{1,3}\)  will be worked out by computing their    corresponding graded ring structures  in terms of generators and relations as well as using the  Lie theoretic data of the associated Lie groups.  Then we use  their graded rings and  Hilbert series  to construct some  families of polarized varieties as weighted completed intersections inside them. The graded rings of the constructed polarized varieties are induced from the graded rings of the ambient key varieties \(w\lgr(3,6)\) and \(w\fl_{1,3}\).

       We construct  projective varieties \((X,D)\)  polarized by \(\QQ\)-ample Weil divisor (i.e. \(nD\) is a Cartier divisor)   \(D \)  with finitely generated graded ring \begin{displaymath}
R(X,D)=\bigoplus_{n \geq 0}H^0(X,nD).
\end{displaymath}
The embedding \[i : X = \Proj R(X,D) \into \PP[w_0, \cdots , w_n]\]  is provided by surjective morphism \begin{displaymath}
\CC[x_0, \cdots , x_n]\onto R(X,D)
\end{displaymath}
from a free graded ring \(S=\CC[x_0, . . . , x_n]\).  The graded ring \(S\) is  generated by variables \(x_i\) of weights \(w_i\). The  divisorial sheaf \(\Oh_X(D)\) of \(X\) is  isomorphic to \(\Oh_X(1) = i^*\Oh_{\PP}(1). \) 

In the past, such examples  have been computed  in~\cite{fletcher,ABR}, where the  cases of codimension at 
most 3 are discussed. The codimension 4 case was initially  discussed in~\cite{altinok} and studied recently more rigorously in~\cite{tom-jerry-1}; see  Corti and Reid \cite{wg} 
for examples in codimension 5 and \cite{qs,qs-ahep}   for  examples in  codimension 6 and 8.

  We construct   examples in codimensions 7 and 9 by taking  quasilinear sections of weighted  flag varieties $(w\Sigma, \Oh_{w\Sigma}(1))$ embedded in weighted projective space \(w\PP V_\la\) 
by their natural Pl\"ucker-type embeddings.  We compute the Hilbert series of a given weighted flag variety to find the canonical divisor class of \(w\Si\). Then we take  quasilinear sections (general hypersurfaces of the appropriate degree in weighted projective space) of \(w\Si\) or of  projective cone(s)  over it to get a  variety with the desired canonical or anticanonical class. We need the defining equations of flag varieties to understand the nature of singularities. The defining ideals of flag varieties are given in~\cite[Sec 1]{rudakov}. We work out the equations by following the algorithmic approach of DeGraaf \cite{degraaf}, coded in computer algebra system  GAP4~\cite[Appendix A]{qs}.  

 We fix the notations and give the necessary definitions in Section~\ref{notationsec}. We also provide 
a quick review of the weighted flag varieties, the  formula for  their Hilbert series $P_{w\Sigma}$, the   defining ideals of flag varieties, and the 
general method of  constructing  families of polarized varieties   as quasi-linear 
sections of \(w\Si\) . In Section \ref{lg36sec}, we study weighted flag varieties associated 
to the symplectic  Lie group  \(\Sp(6,\CC)\), leading to the codimension seven varieties.  The case of
weighted partial flag variety \(w\fl_{1,3}\) in codimension nine is discussed in Section \ref{a3sec}. In Section \ref{classification}, we give the list of all possible distinct flag varieties, which have embeddings   in codimension  \(c\), for  \(4\le c\le 10\).

\subsection*{Acknowledgements}
I am grateful to  Bal\' azs Szendr\H oi  for  helpful discussions.  This research was partially supported  by a grant of  the Higher Education Commission (HEC) of Pakistan. I am also thankful to the International Centre of Theoretical Physics, Italy for their hospitality during the write up of this article.
\section{Definitions and  conventions }\label{notationsec}

We work over a field \(\CC\) of complex numbers. A   pair \((X,D)\), where~\(X\) is a normal projective algebraic variety and~\(D\) is an 
ample $\QQ$-Cartier Weil divisor on $X$ is called polarized variety. All our varieties  appear as projective subvarieties of some weighted projective space. 

 The standard notation
$\PP[w_0,w_1,\cdots,w_n]$  denotes the  weighted projective space; sometimes we will write $w\PP$
if no confusion can arise. The weighted projective space \(\PP^n[ w_i]\) is called well-formed, 
if no \(n-1\) of \(w_0,\cdots,w_n\) have a common factor. 

A polarized variety \(X\subset \PP^n[w_i]\) of codimension \(c\)  is called well-formed, 
if \(\PP^n[w_i] \) is well-formed and \(X\) does not contain a codimension \(c+1\) singular 
stratum of \( \PP[w_i]\). 
\(X\) is called quasi-smooth if the affine cone \(\widetilde X \subset \AAA^{n+1}\) of \(X\) is smooth outside its vertex \(\uline{0}\). If \(X\) is quasi-smooth, then it will only have quotient singularities induced by the  
singularities of \(\PP[w_i]\). We assume that  polarization is provided by the restriction of the tautological ample divisor 
$\Oh_{\PP}(1)$. 

The Hilbert series of a polarized projective variety \((X,D)\) is 
\[ P_{(X,D)}(t)=\sum_{n \geq 0}\dim H^0(X,nD) t^{n}.
\]
We will sometimes write $P_X(t)$ if no confusion can arise. Appropriate Riemann--Roch 
formulas, together with vanishing, can be used to compute $h^0(X,nD)=\dim H^0(X, nD)$ 
in favourable cases. 

A polarized  Calabi--Yau 3-fold \((X,D)\) is a Gorenstein, normal, projective 
three dimensional algebraic variety with \(K_{X} \sim 0\)  and 
\(H^{1}(X,\Oh_X)=H^2(X,\Oh_X)=0\). We allow \((X,D)\) to have at worst canonical quotient 
singularities, consisting of points and curves on \(X\).

\label{sec:flag}
Let  \(G \)  be a reductive Lie group  with its Lie algebra being denoted by   $\mathfrak g$. 
We denote by \(B\) a Borel subgroup of \(G\), by \(P\)  a parabolic subgroup of \(G\), 
and by \(T\) a maximal torus, such that \(T\subset B \subseteq P \subset G\) and \(\got\subset \gob\subseteq\gop \subset \gog\) denote the corresponding Lie algebras inclusions. 

Let $\Lam_W={\rm Hom}(T, \CC^*)$ be the weight lattice and   $V_\la$ denote the $G$-representation with highest
weight \(\la\) where \(\nabla(V)\) represents the set of weights of the representation \(V_{ \la}\). The quotient \(\Si=G/P\) of the Lie group \(G\) by a parabolic subgroup \(P\)  is 
called (generalized or  partial)  flag variety; if \(P=B\) then we call \(\Si\) a complete flag variety.
 The flag  varieties \(\Si=G/P_\la\) are projective subvarieties of  \(\PP V_\la\), where    
$P = P_\la$ is the parabolic corresponding to the set of simple roots \(\nabla_0\) of $G$ orthogonal to the weight vector \(\la\).

 For example, if  \(G=\SL(n,\CC)\), then $B$ is the subgroup of upper triangular matrices. If \(P\) is the parabolic subgroup  containing \(B\),  then  the quotient \[\Si= \SL(n,\CC)/P=\lbrace 0 \subset V_i \subset V_m \subset \cdots \subset V_n \rbrace\] is  a (partial or generalized) flag variety, with \(P\) being the stabilizer of a fixed partial flag with dimension vector \(I=(i,m,\cdots,n)\). Whereas for the dimension vector  \(I=(1,2,\cdots,n),\) we have \(P=B\) and the quotient \(\Si=G/B\) is a complete flag variety.

\subsection*{Weighted flag varieties}\label{hssec}
Let $\Lam_W^* = \Hom(\CC^*,\ZZ)$ be the lattice of one-parameter subgroups of \(G\).  Choose \(\mu \in \Lam_W^* \) and an integer \(u \in\ZZ\) such 
that \[
\left<w\lam,\mu\right>+u >0
\label {weights} 
\]
for all  elements \(w\) of the Weyl group \(W\) of the Lie group \(G\), where $\langle ,\rangle$ denotes the perfect 
pairing between $\Lam_W$ and $\Lam_W^*$.
We recall the definition of weighted flag variety by Grojnowski and Corti--Reid .
\dfn \cite{wg}  
Let \(\Si\) be a flag variety.  
Take the affine cone $ \widetilde{\Sigma} \subset  \widetilde{V_{\lambda}}$ of the 
embedding $\Sigma \hookrightarrow \PP V_{\lambda}$,
then the quotient of    $V_{\lambda}\backslash\{0\}$ by the $ \CC^*$-action given by
\[ 
(\varepsilon \in \CC^*) \mapsto ( v \mapsto \varepsilon^u(\mu(\varepsilon)\circ v))
\] is called weighted flag variety. We also use the term key variety or Format for \(w\Si\). We denote
 this variety by \(w\Si(\mu,u)\) or in short $w\Sigma$, if no confusion can arise. 
 
The Hilbert series of a weighted flag variety \[P_{w\Si}(t)=\sum_{m\ge 0} \dim \left(H^0(w\Si,\Oh(mD))\right)t^m\] can be computed by using the following theorem.
\begin{theorem}\cite[Thm. 3.1]{qs} The Hilbert series of the weighted flag variety 
$\left(w\Sigma(\mu,u),D\right)$ has the following closed form.
\begin{equation}
P_{w\Si}(t)=\dfrac{\displaystyle\sum_{w\in W}(-1)^w \dfrac{t^{\left<w\rho, \mu\right>}}{(1- t^{\left<w\lambda,\mu\right>+u})}}{\displaystyle\sum_{w\in W}(-1)^w t^{\left<w\rho, \mu\right>}}.
\label{whhs}
\end{equation}
Here \(\rho\) is the Weyl vector, half the sum of the positive roots of \(G\), and  $(-1)^w=1  \; \mbox{or} -1$ depending on whether $w$ consists of an even or odd number of simple reflections in the Weyl group $W$, and \(D=\Oh_{w\Si}(1)\) under the embedding \(w\Si \subset w\PP V_{\lam}\).
\end{theorem}

%  Note that by the Weyl denominator identity, the denominator of our expression for 
% the Hilbert series has two equivalent forms 
% \begin{equation}
% \sum_{w\in W}(-1)^w t^{<w\rho,\mu>}=t^{<\rho,\mu>} \prod_{\al \in \nabla_+}(1-t^{<-\al,\mu>}).
% \label{hdidentity}
% \end{equation}
\begin{rmk}  By the standard 
Hilbert--Serre theorem \cite[Theorem 11.1]{atiyah}, 
the Hilbert series of the weighted flag variety has a reduced expression
\begin{equation}
P_{w\Si}(t)=\dfrac{N(t)}{\displaystyle \prod_{\la_i \in \nabla{(V_\la)}}(1-t^{<\la_i,\mu>+u})}.
\label{reducedhs}
\end{equation}
The polynomial \(N(t)\) is called the Hilbert numerator of the Hilbert series and contains some information about the free resolution of the graded ring \(R\left(w\Si,\Oh (mD)\right)\).

\end{rmk}

The flag variety \(\Si=G/P \into \PP V_{\lam}\) is defined by an ideal  \(I =\left<Q\right>\)  
of quadratic equations \cite[2.1]{rudakov}. The  second symmetric power of the contragradient representation \(V^*_{\lam}\) has a  decomposition   
\begin{displaymath}
Z=V_{2\nu}\oplus V_1 \oplus\cdots \oplus V_n
\end{displaymath}  
into irreducible direct summands as a \(G\)-representation, with \(\nu\) being the highest weight of the 
representation \(V^*_{\lam}\). The generators of the linear subspace  \(Q \subset Z=S^2 V^*_{\lam}  \)
 consisting of all the summands except~\(V_{2\nu}\), gives  the defining equations of the flag variety.   
The equations of $w\Sigma$ can be readily computed from this information using computer algebra~\cite{qs}.

\subsection*{Constructing polarized varieties}
We start by constructing  some key variety \(w\Si \) (weighted flag variety) embedded into some weighted projective  space.  Then we find the Hilbert series of the given variety, which gives us some information about the graded ring \(R(w\Si,D)\).  Under suitable  conditions, we can compute   the canonical divisor class of  \(w \Si\). Then we take  quasilinear sections (general hypersurfaces of appropriate degree in weighted projective space) of \(w\Si\) or of  projective cones  over it, to get a  variety with the desired canonical or anticanonical class. Then we study  different aspects, such as  singularities, well-formedness and quasi-smoothness of the resulting variety to establish the existence of an appropriate model of the variety. At the end, 
using the orbifold 
Riemann--Roch formula of~\cite[Section 3]{anita}, we  compute the invariants of our polarized variety \((X,D)\)    from the first few values of $h^0(nD)$, and verify 
that the same Hilbert series can be recovered.
 More  details can be found in \cite{qs}. 
\section{Weighted Lagrangian Grassmannian $ w$LGr(3,6)   varieties} \label{lg36sec}
 \subsection{Generalities}The { symplectic Lie group} $\Sp(6,\CC)$ is the group of automorphisms $A$ of $\CC^6$ preserving a  nondegenerate, skew-symmetric, bilinear form called  the   symplectic form $J$, that is $$J(Av,Aw)=J(v,w) \; \textrm{for all}\; v,w \in V.$$ By using a change of basis, \(J\) has the normal  form:  \[J=\left(\begin{array}{cc}0& I_n\\-I_n&0 \end{array}\right).\] The symplectic group has the embedding in $\GL(6,\CC)$:\[\Sp(6,\CC)=\left\{M\in\GL(6,\CC):M^{t}J M=J\right\}.\]
Let $G$  be the symplectic Lie group  $\Sp(6,\CC)$ with maximal torus $T$.  The weight lattice \(\Lam_W\) of \(\Sp(6,\CC)\) is  a rank 3 lattice $\Lam_W=\left<e_1,e_2,e_3\right>.$ The simple roots in the weight lattice of the Lie algebra $\gosp_6$  are   $$\al_1=e_1-e_2,\;\al_2=e_2-e_3,\mbox{ and }\al_3=2e_3.$$  The dominant fundamental weights of $ \Sp(6,\CC) $ are given by $ \om_i= e_1+\cdots+e_i$, for $1 \leq i \leq 3.$ The sum of the fundamental weights, also known as Weyl vector, is given by  $\rho=3e_1+2e_2+e_3.$

Consider  a \(G\)-representation  \(V=\bigwedge^3\CC^6\); the  third exterior power of the  standard representation \(\CC^6\), which is 20 dimensional. Then we have a natural contraction map \(\kappa:\bigwedge^3\CC^6\to \CC^6, \) obtained by contracting with elements of $\bigwedge^2\left( \CC^6\right)^*$. The map $\kappa$ decomposing     \(V\) into its  summands  \(V=V_\la\oplus V_1\), the kernel of  \( \kappa\) is  an irreducible representation $V_\la$ of \(G\) with  highest weight $\lam=\om_3=e_1+e_2+e_3$. Then  the Weyl dimension formula  tells us that  $V_{\lam}$ is 14 dimensional. All the fourteen weights of \(V_{\la}\) appear  with multiplicity one. If   \(\nabp=\left\lbrace  \al_1,\al_2,\al_1+\al_2\right\rbrace,  \) then 
  the corresponding   parabolic subalgebra   \[\gop_{\lam}= \bigoplus\left(\got \bigoplus_{\al \in \nabla_+}\gog_{\al}\bigoplus_{\al \in {\nabp}} \gog_{-\al}\right)\]  is 3+9+3 = 15 dimensional. Thus  the corresponding flag variety     \(\Si=G/P_{\lam}\) is six dimensional:
   \[\dim(\Si)=\dim(\gosp_6)-\dim(\gop_{\lam})=21-15=6.\]   This flag variety is known as the Lagrangian Grassmannian  $\lgr(3,6) $ of Lagrangian subspaces in $\CC^6$ (\cite{harris}) and also as the {\it symplectic Grassmannian} ( \cite{mukai6}).  Obviously we have a codimension seven    embedding \(\lgr(3,6) \into \PP^{13} V_\la\).

\subsection{The weighted flag variety}\label{lgr36-const}The weighted flag variety $ w\lgr(3,6) $ can be constructed by letting   $ f_1,f_2,f_3 $ be the basis of the  dual lattice $ \Lam_W^* $, dual to \(e_1,e_2  \mbox{ and } e_3.\) Then for $ \mu= a_1f_1+a_2f_2+a_3f_3 \in \Lam_W^*$ and $ u \in \ZZ $, to  get the weighted version of $ \lgr(3,6) $:
\[ w\Si(\mu,u)=w\lgr(3,6) \into w\PP^{13}. \]  
  The set of weights on weighted projective space $w\PP^{13}$ is  $\{<\la_i,\mu>+u\}$, where $\la_i$ are the weights of the representation \(V_\la\).  We  denote an element \(\mu \) of the  dual lattice \(\Lam_W^*\) by a vector of integers, i.e. $\mu=(a_1,a_2,a_3).$

\subsection{Hilbert Series of $w$LGr(3,6)} 

\thm Let \(S_2\) be the symmetric group on 2 elements and \(\si(a)=a\) if \(\si\) is  even and $\si(a)=-a$ if \(\si\) is  odd permutation. Then  the Hilbert series of the \(w\lgr(3,6) \)  has the compact form
\begin{equation}\label{hslg36} 
P_{w\lgr(3,6)}(t)=\dfrac{1-P_1(t)\left( t^{2u}-t^{9u}\right)+P_2(t)\left(t^{3u}-t^{7u}\right)-P_3(t)(t^{4u}-t^{6u}) -t^{10u}}{\displaystyle\prod_{\la_i \in\nabla(V_\la)} \left( 1-t^{<\la_i,\mu>+u} \right)},
\end{equation}
where 
 \[\begin{array}{ll}
P_1(t)&=\displaystyle\sum_{1 \leq(i,j)\leq 3}t^{a_i-a_j}+\sum_{\sigma\in S_2}\displaystyle\sum_{1 \leq i \leq j \leq3} t^{\si(a_i+a_j)}\end{array},
\]
 \[ \begin{array}{ll}
 P_2(t)&= \displaystyle\sum_{\sigma\in S_2}\left( \sum_{ 1 \leq i < j  \leq 3}\left(t^{ \si (2a_i+a_j)}+t^{\si (2a_i-a_j)}+t^{\si (a_i+2a_j)}+ t^{\si (a_i-2a_j)}\right)\right.\\
&\left.+2\left(t^{\si (a_1+a_2+a_3)}+t^{\si (a_1+a_2-a_3)}+t^{\si (a_1-a_2-a_3)}+
t^{\si (a_1-a_2+a_3)}\right)
+4  \displaystyle\sum_{i=1}^3 t^{\si (a_i)}  \right)
\end{array} ,\]and 

  \[\begin{array}{ll}P_3(t)=&\displaystyle \sum_{\si \in S_2}\left( \sum_{i=1}^3t^{2\si (a_i)}+3\sum_{1\leq i<j \leq3}^3\left(t^{\si (a_i-a_j)}+t^{\si (a_i+a_j)}\right)\right)\\&+\displaystyle \sum_{\si \in S_2} \left(t^{\si(a_1+ 2a_2)}+t^{\si(2a_1+ a_2)}+t^{\si(a_1- 2a_2)}+t^{\si(2a_1- a_2)}\right)
\sum_{\si \in S_2}t^{\si(a_3)})
\\& +\displaystyle \sum_{\si \in S_2}\left( t^{\si(a_1+a_2)}+t^{\si(a_1-a_2)}\right)\sum_{\si \in S_2}t^{2\si(a_3)} +4\end{array}.\] Moreover, if 
\(w\lgr(3,6)\) is  well-formed  then the  canonical line bundle \(K_{w\lgr(3,6)}=\Oh_{w\lgr(3,6)}(-4u).\)
\proof The Weyl group \(W\) of the symplectic Lie algebra \(\gosp_6\) is the semidirect product of the symmetric group  \(S_3\) with 3 copies of \(\ZZ/2\ZZ\); \[W= S_3 \ltimes (\ZZ/2\ZZ)^3.\]  Using representation theory  we work out  that eight of the fourteen weights appear as the image of \(\lam \) under the action of  the Weyl group \(W\). Each of the eight weights appears with multiplicity six in the orbit of \(\la\) under the  \(W\)-action, accounting for the 48 elements of the Weyl group \(W\). After  performing  some simplifications we get   the following form of the Hilbert series of  $ w\lgr(3,6).$
\begin{equation}\label{lgr36ie}P_{w\lgr(3,6)}(t)= \dfrac{1+\displaystyle \sum_{\si \in S_2} \sum_{1 \leq i \leq 3}t^{\si (a_i)} \left( t^u-t^{3u} \right)-t^{4u}}{\displaystyle\prod_{\la_i \in W\la} (1-t^{<\la_i,\mu>+u})}
,\end{equation} where \(W\la\) is the orbit of \(\la\) under \(W\)-action. 
 As the weights \(a_i+u,-a_i+u\) with \(1\le i\le 3\) do not appear  in the orbit of \(\la\) under the \(W\)-action, we multiply and divide the expression \eqref{lgr36ie} with      \[ P_L(t)=\prod_{1 \leq i \leq3} \left(1-t^{a_i+u}\right)\left(1-t^{-a_i+u}\right ),\] to get the full expression for the Hilbert series of \(w\lgr(3,6)\).
 After performing some  simplifications we get the required compact form \eqref{hslg36}. As the adjunction number of \(w\lgr(3,6) \) is $ 10u $ and the sum of weights on $ w\PP V_{\lam} $  is  $ 14u $;  the canonical divisor is given by    \[ K_{w\lgr(3,6)}= \Oh_{w\lgr(3,6)}(10u-14u)=\Oh_{w\lgr(3,6)}(-4u).\square\]
\begin{comment}\begin{rmk}  From the Hilbert series and the formation of weights, it is clear that we can choose the parameter \(\mu\) to be half integer, if we choose the each of the coefficients of \(\mu\) to be half integer as well.    
\end{rmk} \end{comment}
\subsection{Examples} 
 
\example Consider the Hilbert series of the straight flag variety $ \lgr(3,6) $, which corresponds to $ u=1 $ and $ \mu= \underline{0}. $ Then we have \[P_{\lgr(3,6)}(t)= \dfrac{1-21t^2+64t^3-70t^4+70t^6-64t^7+21t^8-t^{10}}{ \left ( 1-t \right)^{14}}.\] Since $ \lgr(3,6) $ is a 6-dimensional well-formed and smooth variety, we can compute the  canonical bundle \(K_{\lgr(3,6)}\) to be $ \Oh(-4) $. Let \(H_1, H_2 \text{ and } H_3\) be three general hyperplanes of \(\PP^{13}\) then we get a  three dimensional variety  \[V=\lgr(3,6) \cap H_1 \cap H_2 \cap H_3 \subset \PP^{10}\] with \(K_{V}=\Oh_V(-1). \) Thus \(V\) is  a Fano 3-fold of genus 9, anti-canonically polarized by \(K_V\)  where     $ (-K_V)^3=16 $. This variety was constructed by Mukai by using the vector bundle method in \cite{mukai}.  

\begin{rmk}We searched for more families  of terminal $ \QQ $-Fano threefolds but  we did not manage to find a new family in codimension 7. A list of 303 such  $\QQ$-Fano 3-folds 
with terminal singularities can be found  on Gavin Brown's graded ring data base page \cite{grdb}.
 \end{rmk}
\begin{example}\label{gr26exp1} Consider the following initial data
 \begin{itemize}
\item Input: $ \mu=(1,0,0) $, $ u=2 $
\item Variety and weights: $ w\lgr(3,6)\subset \PP^{13}[1^5,2^4,3^5]  $, with weights assigned to the variables $x_i$ in the equations given in appendix \ref{app:lgr36} on Page \pageref{eqlgr36}.
\[
 \renewcommand{\arraystretch}{1.5}
 \begin{array}{ccccccccccccccc}
{\rm Variable} & x_1   & x_2   & x_3   & x_4&x_5&x_6&x_7&x_8&x_9&x_{10}&x_{11}&x_{12}&x_{13}&x_{14}\\
{\rm Weight}  &  3 &  3&           3&    3  &  2 &  3 &      2&    2&     1 &   2&    1&   1&    1&    1
 \end{array}
\]
\item Canonical class: $ K_{w\lgr(3,6)}= \Oh(-8) $
\item Hilbert Numerator: $1-t^2-4 t^3-7 t^4+12 t^5+\cdots-12 t^{15}+7 t^{16}+4 t^{17}+t^{18}-t^{20}$
\end{itemize} 
We take a  threefold quasilinear section \[ X=w\lgr(3,6) \cap  (3)^2\cap (2) \subset \PP^{10}[1^5,2^3,3^3] ,\] then \[ K_X= \Oh_X\left(-8+(2\times3+2 )\right) \sim 0.\]
By using the equations, we can check that  $ X $  does not contain any codimension 8 singular strata. Therefore $ X $ is a well-formed polarized weighted projective threefold. We work out  the singularities  of $ X $ induced by the weights of $w\PP^{10}$. To make a choice of coordinates we suppose  
\[\{x_3=f_3 (x_i), \; x_4=g_3(x_i) \;,  x_5=f_2(x_i)\}.\] Here $ f_3,g_3 $ and $ f_2 $ are general  homogeneous equations of degree $3, 3$ and 2 respectively  in the rest of the variables. Now we check the singularities of $ X $ on each of the singular strata and find the local transverse structure of  these singularities.\\
\textbf{$1/3$ singularities} The defining locus $ C $ of this singular stratum is defined by taking $$C= X \cap \left\lbrace  x_7=x_8=x_9=x_{10}=x_{11}=x_{12}=x_{13}=x_{14}= 0  \right\rbrace .$$ Therefore $$ C:= \left \lbrace \frac{f_3^2}{4}-g_3x_2+x_1x_6 =0 \right \rbrace  \subset \PP^2[x_1,x_2,x_6].$$
This defines a  quadratic  curve in \(\PP^2\). On each point of the curve one of the \(x_i \neq 0\). If we consider  \(x_1 \neq 0,  \) then by using the implicit function theorem we can eliminate the variables \(x_{6},x_8,x_{10},x_{11},x_{12},x_{13}\) and \(x_{14}\) from equation number (A1), (A2), (A4), (A3), (A5), (A6) and (A7) respectively appearing in appendix \ref{app:lgr36} on page \pageref{eqlgr36}. Therefore \(x_7\) and \(x_9\) are local variables near this point and the group \(\mu_3\) acts by \[ \ep : x_7,x_9 \mapsto \ep^2 x_7, \ep x_9 .\] For \(x_2\neq 0\), we can eliminate the variables \(x_7,x_9,x_{10},x_{12},x_{13}\) and \(x_{14}\) from equation number (A2), (A3), (A8), (A9), (A10) and (A11) respectively. If we suppose that \(g_3\) contains the monomial \(x_6\), then we eliminate \(x_6\) from equation (1) as well. Therefore \(x_{8}\) and \(x_{11}\) are the local variables on this point of the curve and the group \(\mu_3\) acts by \[\ep:x_8,x_{11}\mapsto \ep^2x_8,\ep x_{11}.\] Similarly, near the point \(x_6 \neq 0\) we work out that \(x_{10}\) and \(x_{14}\) are local variables and the group \(\mu_3\) acts by \[\ep:x_{10},x_{14}\mapsto \ep^2x_{10},\ep x_{14}.\] 
Therefore \(C\) is a rational curve of singularities of type \(\dfrac{1}{3}(1,2)\).\\
\textbf{\(1/2\) singularities:} This singular locus is defined by restricting \(X\) to the locus defined by \(x_7,x_8\) and \(x_{10}\).  The equation (A6)
 contains the monomial \(x_{7}^2\) and no term involving \(x_8\) and \(x_{10}\), the equation (A11) contains a monomial \(x_{8}^2\) and no term involving \(x_7\) and \(x_{10}\) and equation (A19) contains a monomial \(x_{10}^2\) and no term involving \(x_7\) and \(x_8\). Therefore, we conclude that \(X\) has no singularities along the singular stratum defined by weight 2 variables.\\
\textbf{Quasismooth:} As we have shown, on \(\dfrac13\)  and \(\dfrac12\) strata that  \(X\) is locally a threefold by using the implicit function theorem. Similarly we can show that on the rest of the strata that \(X\) is locally a threefold.  \\
Thus  $(X,D)$ is a polarized Calabi--Yau threefold with  a rational  curve \(C\) of singularities of  type \(A_2\).  The rest of the invariants of this family, computed  by using the orbifold  Riemann--Roch formula of~\cite[Section 3]{anita} for Calabi--Yau threefolds are given as follows.
\begin {itemize}
 \item $D^3= \dfrac{64}{9}$, $D.c_2(X)=48$, 
 $\deg D\rvert_{C}= \dfrac{2}{3}$, $\gamma_C=-2$
\end{itemize}
\end{example}

\section{Weighted  $A_3$ type partial flag  variety}\label{a3sec}
\subsection{Construction of weighted ${\rm FL}_{1,3}$ variety} \label{basicsa3}Let \(G\) be the reductive Lie group  \(\GL(4,\CC)\)  with  corresponding Lie algebra \(\gog=\gogl_4\),  known as    \(A_3\)-type Lie algebra. The rank of the maximal abelian subalgebra \(\got \) is 4, which corresponds to the maximal torus \(T\) inside \(G\).  The weight lattice of the Lie algebra \(\gogl_4\) is  a rank 4 lattice \(\Lam_W=\left<e_1,e_2,e_3,e_4\right> \).  The set of simple roots of the  root system of the corresponding simple part of \( \gog\) is\[\al_i=e_i-e_{i+1}\text{ for }1\leq i \leq 3.\] The Weyl group \(W\) of \(G\) is \(S_4\), the symmetric group on 4 letters, of order 24. The Weyl vector can be taken to be \(\rho=3e_1+2e_2+e_3.\)

Consider   two 4 dimensional representations of \(G\); \(V_{1}=\CC^4\) which is the standard representation of \(G\) and the representation  \(V_2=\bigwedge^3\CC^4\). They are irreducible representations of  \(G\) with  highest weights \(\la_1=e_1\) and \(\la_2=e_1+e_2+e_3.\) Then the representation \(V=V_1\otimes V_2\) is  a 16 dimensional representation of \(G\)  which is not irreducible.  Consider the wedge product map \[\mathcal{S}:\CC^4 \otimes \bigwedge^3 \CC^4 \longrightarrow \bigwedge^4 \CC^4.\]Then   the    kernel \(\kappa(S)\) of the map, is the  irreducible highest weight representation  of \(G\) with    highest weight   \[\la=\la_1+\la_2=2e_1+e_2+e_3.\] By using  the Weyl dimension formula  we can show that \(V_{ \la}\) is 15 dimensional.  Twelve of the weights of \(V_\la\)  appear with multiplicity one and one with multiplicity three. 
The dimension of  the Lie algebra \(\gogl_4\) is 16. The only simple root orthogonal to the highest weight \(\lam\)  under in the weight lattice is   \(\al_2\).  Therefore, the  parabolic subalgebra  \[\gop_{\lam}= \bigoplus\left(\got  \bigoplus_{\al \in \nabla_+}\gog_{\al}\bigoplus \gog_{-\al_2}\right),\]is 4+6+1=11 dimensional. Hence  the corresponding   flag variety \(\Si=G/P_{\lam}\) is  five dimensional and we get  a  codimension 9 embedding \(\Si^5 \into \PP^{14}[V_{\lam}].\)
In the notation of Section \ref{sec:flag}, we have a dimension vector \(I=(1,3,4)\). The reductive Lie group \(\GL(4,\CC)\) parameterizes  the flags of type \[\{0\subset V_1\subset V_3\subset V_4\}.\] 
We will denote this partial flag variety by  \(\fl_{1,3}\). 
 To obtain the  weighted version of \(\fl_{1,3}\), we let \(\Lam_W^*=\left<f_1,f_2,f_3,f_4\right>\) to be the dual lattice of the weight lattice. Then for any   \[\mu= \displaystyle\sum_{i=1}^4 a_i f_i  \in \Lam_W^* \text{ and } u \in \ZZ\]   we get the embedding \begin{equation}\label{fl13-embedding}w\fl_{1,3}(\mu,u)\into w\PP V_{\lam}[\left<\la_i,\mu\right>+u ],\end{equation} where \(\la_i\) are the weights of the representation \(V_{\lam}\) understood with multiplicities. Following the convention of the Section \ref{lgr36-const},   element \(\mu \) of the  dual lattice \(\Lam_W^*\) is represented by  $\mu=(a_1,a_2,a_3,a_4). $  Here we can choose all \(a_is\) to be the half or quarter integers to get the embedding \eqref{fl13-embedding}  but we can get the embedding with same set of weights  by  changing the value  of \(u\), so we only take them to be integers. 
\subsection{Hilbert series of weighted ${\rm FL}_{1,3}$ }  
\thm  Let \(s\) be the sum of the integers in \(\mu=(a_1,a_2,a_3,a_4). \) Then the  Hilbert series of the  weighted flag variety \(w\fl_{1,3}\)  has the following compact form.
\begin{equation}\label{hsa3}P_{w\fl_{1,3}}(t)=\dfrac{1+\displaystyle\sum_{k=1}^4(-1)^kP_k(t)t^{ks+(k+1)u}+\sum_{k=5}^8(-1)^kP_k(t)t^{(k+1)s+(k+2)u}+t^{12(s+u)}}{\displaystyle\prod_{w_i \in\nabla( V_{\lam})}\left(1-t^{<w_i,\mu>+u}\right)},\end{equation} \rm{where} 
\[\begin{array}{l}P_1(t)=\displaystyle\sum_{1 \leq i <j \leq 4}t^{2(a_i+a_j)}+2\sum_{1 \leq(i,j)\leq 4}\left (t^{s+a_i-a_j}-t^{s}\right),\end{array}\] 
\[\begin{array}{ll}P_2(t)&= 4\displaystyle\sum_{1 \leq i < j \leq 4}t^{2(a_i+a_j)}+8\sum_{1 \leq (i,j)\leq 4}t^{(a_i-a_j)+s}+\\&\displaystyle\sum_{1 \leq (i,j,i \neq j) \leq 4}\left(t^{2s-(3a_i+a_j)}+ t^{3a_i+a_j}\right)-16t^{s} \end{array},\]
\[\begin{array}{ll}P_3(t)&= 6\displaystyle\sum_{1 \leq i < j \leq 4}t^{2(a_i+a_j)}+14\sum_{1 \leq (i,j)\leq 4}t^{(a_i-a_j)+s}\\&+\displaystyle\sum_{1 \leq (i,j,i \neq j) \leq 4}\left(3t^{2s-(3a_i+a_j)}+3 t^{3a_i+a_j}+t^{2(a_i-a_j)}\right)-29t^{s} \end{array},\]  

\[\begin{array}{ll}P_4(t)&= 4\displaystyle\sum_{1 \leq i < j \leq 4}t^{2(a_i+a_j)}+12\sum_{1 \leq (i,j)\leq 4}t^{(a_i-a_j)+s}\\&+\displaystyle\sum_{1 \leq (i,j,i \neq j) \leq 4}\left(3t^{2s-(3a_i+a_j)}+3 t^{3a_i+a_j}+2t^{2(a_i-a_j)+s}\right)-24t^{s} \end{array}.\] If \(w\fl_{1,3}\) is well-formed then \[K_{w\fl_{1,3}}=\Oh_{w\fl_{1,3}}(-3(s+u)). \]

\proof The Weyl group \(W\) of the symmetric group \(S_{4}\). By using representation theory  we work out  that twelve  of the fifteen  weights are in the orbit of  \(\lam \) under the action of  the Weyl group. To compute the Hilbert series of \(w\fl_{1,3}\), we evaluate the  expression \eqref{whhs} for   \(W,\rho, \mu\) and \(\lam\) as described in  Section \ref{basicsa3}. After  some simplification we get the following form of the Hilbert series of the weighted  flag variety \(w\fl_{1,3}\). \begin{equation}\label{a3hs} P_{w\fl_{1,3}}(t)= \dfrac{1+3 t^{s+u}-\left(\displaystyle\sum_{1 \leq i <j \leq 4}t^{2(a_i+a_j)}+\sum_{1 \leq(i,j,i\neq j)\leq 4}2t^{s+a_i-a_j}\right)t^{s+2u}+\cdots+t^{9(s+u)}}{\displaystyle\prod_{\la_i \in W\la} 1-t^{<\la_i,\mu>+u}},\end{equation}
where \(s=\displaystyle\sum_{i=1}^4a_i, \) and \(W\la\) denote the orbit of \(\la\) under \(W\)-action.   The full expression for the Hilbert series of   \(w\fl_{1,3}\) is obtained by multiplying and dividing \eqref{a3hs} by \((1-t^{s+u})^3\), which represent those  weight spaces of \(\PP^{14}[\left<\la_i,\mu\right>+u]\) which do not lie in the orbit of $ W$-action.   After further simplifying we get the required compact form \eqref{hsa3}  of the  Hilbert series of \(w\fl_{1,3}\). The adjunction number is \(12(s+u)\) and the sum of the weights on \(w\PP V_\la\) is \(15(s+u),\) thus   \(K_{w\fl_{1,3}}= \Oh_{w\fl_{1,3}}(-3(s+u)). \)$\square$
\begin{rmk}Due to the  Gorenstein symmetry of the resolution of the  graded ring corresponding to  \(w\fl_{1,3}\), we have \[P_i(t)=P_{c-i}(t), \text{where \;\;}c=\codim(w\fl_{1,3}),\; i =1,2,3,4.\] The 36 defining quadrics are visible in the Hilbert numerator of the Hilbert series of \(w\fl_{1,3}\) from \(P_1(t).\)
\end{rmk}  
\subsection{Examples } 
\example First we consider the case of the straight partial  flag variety \(\fl_{13}\). We evaluate the expression \eqref{hsa3} for \(\mu=\underline{0}\) and \(u=1\)  to get the  Hilbert series of  the  straight \(A_3\) flag variety \(\fl_{1,3}\); 
 \[P_{\fl_{1,3}}(t)=\dfrac{1-36 t^2+160 t^3-315 t^4+288 t^5-288 t^7+315 t^8-160 t^9+36 t^{10}-t^{12}}{(1-t)^{15}}.\] The canonical divisor  class of  ~\(\fl_{1,3} \)  can be read off from the Hilbert series, which is   \(K_{\fl_{1,3}}= \Oh_{\fl_{1,3}}(-3).\) Then the  intersection of \(\fl_{1,3}\) with two general hyperplanes \[V= \fl_{1,3} \cap H_1 \cap H_2 \subset \PP^{12}\] is a  non-prime Fano 3-fold, polarized by its anti-canonical class \(-K_V\).  The degree of the embedding is given by \(\left(-K_V\right)^3=20\) and by using \((-K_V)^3=2g-2\) it is evident that the genus of \(V\) is 11.  
This variety is listed in \cite{mukaim} as non-prime  Fano threefold. The non-primeness also follows from the fact that $\Pic(\Si)=\ZZ^2,$ by using standard theory \cite[Sec 6.3]{BE}. Our description of this variety as a linear section  of a partial flag variety  does not seem to appear in any of  Mukai's articles. 

\begin{rmk}We searched for more families  of terminal $ \QQ $-Fano threefolds but  like previous attempts  we did not manage to find a new family in codimension 9. A list of all possible 93    $\QQ$-Fano 3-folds 
with terminal singularities can be found  on Gavin Brown's graded ring data base page \cite{grdb}.\end{rmk}
\begin{example} Consider the following initial data
\begin{itemize}
\item Input: $ \mu=(0,0,1,1) $, $ u=0 $
\item Variety and weights: $ w\fl_{1,3}\subset \PP^{14}[1^4,2^7,3^4]$, with weights assigned to the variables $x_i$ in the equations given in appendix \ref{eqa3cd9}.
\[
 \renewcommand{\arraystretch}{1.5}
 \begin{array}{cccccccccccccccc}
{\rm Variable} & x_1& x_2& x_3 & x_4&x_5&x_6&x_7&x_8&x_9&x_{10}&x_{11}&x_{12}&x_{13}&x_{14} & x_{15}\\
{\rm Weight}  &  1 & 1& 1& 2 &  2 &  1 &      2 &   2&    2&   3&    2&    2&3&3&3
 \end{array}
\]
\item Canonical class: $ K_{w\fl_{1,3}}= \Oh(-6) $
\item Hilbert Numerator: $ 1-t^2-8 t^3-10 t^4+32 t^5+\cdots-32 t^{19}+10 t^{20}+8 t^{21}+t^{22}-t^{24} $
\end{itemize} 
 We take a projective cone over \(w\fl_{1,3}\), so that we get the embedding \[\Cc w\fl_{1,3} \subset \PP^{15}[1^5,2^7,3^4], \mbox{ and } K_{\Cc w\fl_{1,3}}=\Oh(-7).\]We take a  threefold complete intersection \[ X= \Cc w \fl_{1,3}\cap  (2)^2\cap (3) \subset \PP^{12}[1^5,2^5,3^3] ,\] then \( K_X \sim 0.\)
Then   $ (X,D) $  is a well-formed and  quasismooth Calabi--Yau threefold  with a rational curve of singularities of type \(\dfrac{1}{3}(1,2) \). The rest of the invariants of \((X,D)\) are listed below.
\begin{itemize}  
 \item[$\spadesuit$] $D^3= \dfrac{76}{9}$, $D.c_2(X)=48$, $\deg D\rvert_{C}= \dfrac{2}{3}$, $\gamma_C=10$
\end{itemize}
\end{example}
\begin{example} Initial data
 \begin{itemize}
\item  Input: $ \mu=(0,1,1,1) $, $ u=2 $
\item Variety and weights: $ w\fl_{1,3}\subset \PP^{14}[1^3,2^9,3^3]$, with weights assigned to the variables $x_i$ in the equations given in Section \ref{eqa3cd9}.
\[
 \renewcommand{\arraystretch}{1.5}
 \begin{array}{cccccccccccccccc}
{\rm Variable} & x_1& x_2& x_3 & x_4&x_5&x_6&x_7&x_8&x_9&x_{10}&x_{11}&x_{12}&x_{13}&x_{14} & x_{15}\\
{\rm Weight}  &  1 & 1& 2& 1 &  2 &  2 & 2 & 2&    2&   3&    2&3&2&3&3
 \end{array}
\]
\item Canonical class: $ K_{w\fl_{1,3}}= \Oh(-6) $
\item Hilbert Numerator: $ 1-9 t^3-15 t^4+33 t^5+58 t^6-\cdots-58 t^{18}-33 t^{19}+15 t^{20}+9 t^{21}-t^{24} $
\end{itemize} 
 We take a projective cone over \(w\fl_{1,3}\), so that we get the embedding \[\Cc w\fl_{1,3} \subset \PP^{15}[1^4,2^9,3^3], \mbox{ and } K_{\Cc w\fl_{1,3}}=\Oh(-7).\]Consider  a  threefold quasilinear section \[ X= \Cc w \fl_{1,3}\cap   (2)^2\cap (3) \subset \PP^{12}[1^4,2^7,3^2] ,\] then  the canonical bundle \(K_X \sim 0.\)
Then    $ (X,D) $  is a well-formed and  quasismooth polarized Calabi--Yau threefold  with two rational curves of singularities, \(C\) of type \(A_2\) and \(E\) of type \(A_1\). The rest of the invariants of \((X,D)\) are listed below.
\begin{itemize}  
 \item[$\spadesuit$] $D^3= \dfrac{127}{18}$, $D.c_2(X)=46$, $\deg D\rvert_{C}= \dfrac{1}{3}$, $\gamma_C=8,\;$
 $\deg D\rvert_{E}= 3$, $\gamma_E=1$
\end{itemize}
\end{example}
\section{List of flag varieties in  codimension $4\le c\le 10$ }
%&latex
\label{classification}
Given a  reductive Lie group \(G\) there is one-to-one correspondence between flag varieties $\Si=G/P_\la$ and irreducible highest weight representations $V_\la$.   We  list  all  highest weights $\la$ which leads to the embedding of  $\Si$ in   codimension  $4\leq c \leq 10$.  The list  does  not qualify to be the complete classification of such flag varieties but one can deduce it to be the one, with a bit more effort.   We only list  those $\la$ which lead  to a distinct  embedding  of some flag variety $\Si$ in $\PP V_\la.$ 

We call a flag variety $\Si$ to be  {\it distinct} if   
 \begin{enumerate}
\item it  cannot be recovered as a complete linear section of some  other flag variety $\Si_1=G_1/P_{\la_1}$;
\item it is not isomorphic to some other flag variety $\Si_1=G_1/P_{\la_1}$. \end{enumerate} 
For example the Lagrangian Grassmannian \( \lgr(2,6)\) is a hyperplane  section of the the Grassmannian \(\Gr(2,6)\) in codimension 6. Thus we will only consider  $\Gr(2,6)$ as a variety in codimension 6.
We briefly discuss the possibilities  in each codimension \(c\), for \(4\leq c\leq 10.\)    
\begin{itemize} 
\item \(c=4.\) We have a flag  variety of the Lie group \(\GL(3,\CC)\) with highest weight \(\la=2e_1+e_2\). The corresponding fixed flag parameterised by \(G\)  is \(\{0\subset V_1\subset V_2\subset V_3\}\). Thus we have a complete flag variety \(\fl_{1,2}\).
\item \(c=5.\)  For \(G=\so(10,\CC)\) we have a codimension five  embedding of \[\ogr(5,10)\subset \PP^{15}(S^+)\]for highest weight \(\om_5\), discussed in \cite{wg}. For \(G=\so(9,\CC)\) we have  a codimension five homogeneous variety \(\ogr(4,9)\) corresponding to highest weight \(\om_4\) which  is isomorphic to \(\ogr(5,10).\)
\item \(c=6\). For \(G=\GL(6,\CC)\) and highest weight \(\om_2\) we get a codimension 6 embedding of the Grassmannian  \(\Gr(2,6)\into \PP^{14}\),  discussed in \cite{qs}. For the symplectic Lie group \(\Sp(6,\CC)\) we have a codimension six  embedding of the isotropic   Lagrangian Grassmannian \(\lgr(2,6)\into \PP^{13}\) for the highest weight \(\om_2\), which  is  just a general linear section of \(\Gr(2,6).\)
\item $c=7$. There is only one case which we discussed above in Section \ref{lg36sec}.
\item \(c=8.\) For \(G=G_2\) we have a codimension 8 embedding \(\Si^5\into \PP^{13}\) corresponding to the highest weight \(\om_2\), which is already discussed in \cite{qs}. \item \(c=9.\) The only case  is  discussed in the  Section \ref{a3sec} above. \item \(c=10.\) In this case we  get three distinct homogeneous varieties.\begin{enumerate}
\item Let \(G\) be the  simple Lie group of type \(E_6\) and \(V_\la\) be  the highest weight representation of \(G\) with \(\la=\om_1 \mbox{ or } \om_6\). Then the corresponding homogeneous variety \(\Si\) has a codimension 10 embedding \(\Si \into \PP^{26}V_{\la}\). The  \(F_4\) homogeneous  variety corresponding to the highest weight representation \(V_{\la}\) with \(\la=e_1\), discussed in \cite[Chapter 7]{Q}, is just a general linear section of this \(E_{6}\) variety.
\item For  \(G=\GL(7,\CC)\) we have a representation \( \bigwedge^2\CC^7 \), which leads to the Pl\"ucker embedding of Grassmannian \(\Gr(2,7) \into \PP^{20}; \)  which is a  codimension 10 embedding. \item For  \(G=\GL(6,\CC)\) we have a representation \( \bigwedge^3\CC^6 \), which leads to the Pl\"ucker embedding of the  Grassmannian \(\Gr(3,6) \into \PP^{19}, \) which is also a  codimension 10 embedding. \end{enumerate}  

\rmk Computing  the Hilbert series  and expected syzygies   does suggest that $\Gr(3,6)$ might be a  linear section of $\Gr(2,7) $  as they have exactly the same Hilbert numerator, propagating the same number of generators and syzygies in each degree. In fact it is not the case,  the $H^6(\Gr(2,7))$ is a 2-dimensional and $H^6(\Gr(3,6))$ is a 3-dimensional, so by Lefschetz hyperplane section principle   $\Gr(3,6)$ is not a hyperplane section of $\Gr(2,7)$.        
\end{itemize}
In the following table  we list   all the distinct flag varieties in  codimension \(4\leq c \leq 10\). The last column represents  the number of defining equations of  flag variety $\Si$.  
\begin{table}[h]
\begin{center}\caption{ List of  flag varieties in codimension $4\le c\le$ 10} \label{tab:class}\

\begin{tabular}{|c| c| c| c| c|}  % centered columns (2 columns)
\hline\hline %inserts double horizontal lines
 \(c\)&  Type of \(G\)  & \( \la\)   & Embedding & Number of Eqs\\
 \hline 4&\(\GL(3)\) &\(\om_1+\om_2\)& \(\fl_{1,2}\subset\PP^{7}\)&\(9\)\\
 \hline 5&\(\so(10)\) &\(\om_5\)& \(\ogr(5,10)\subset \PP^{15}\)&\(10\)\\
 \hline 6&\(\GL(6)\) &\(\om_2\)& \(\Gr(2,6)\subset \PP^{14}\)&\(15\)\\
  \hline 7&\(\Sp(6)\) &\(\om_3\)& \(\lgr(3,6)\subset \PP^{13}\)&\(21\)\\
 \hline 8&\(G_2\) &\(\om_2\)& \(\Si\subset \PP^{13}\)&\(28\)\\
 \hline 9&\(\GL(4)\) &\(\om_1+\om_3\)& \(\fl_{1,3}\subset \PP^{14}\)&\(36\)\\  \hline 10&\(E_6\) &\(\om_1\)& \(\Si\subset \PP^{26}\)&27\\
  \hline 10&\(\GL(7)\) &\(\om_2\)& \(\Gr(2,7)\subset \PP^{20}\)&35\\
  \hline 10&\(\GL(6)\) &\(\om_3\)& \(\Gr(3,6)\subset \PP^{19}\)&35\\
\hline
\end{tabular}
\end{center}
\end{table}

\appendix \section{Equations of the Lagrangian Grassmannian LGr(3,6)}\label{app:lgr36}
We  compute the defining   equations of \(\lgr(3,6)\) by the GAP4 code given in Appendix of  \cite{qs}. The code finds the decomposition of the second symmetric power of the dual highest weight representation \(V_\la^*\) into its direct summands as a module over the symplectic Lie algebra \(\gosp_6\).  The second symmetric power of the dual representation \(V_{\lam}^*\) has a decomposition $S^2V_{\lam}^*=\bigoplus (\oplus V_{1} \oplus V_{e})$ into 84 and 21 dimensional representations. The 21 quadrics   generate the 21 dimensional subspace \(V_{e}\) of \(S^2V_{\lam}^*\) and they are listed as follows.  
\begin{eqnarray} \label{eqlgr36}
&&\textstyle  x_1x_6-x_2x_4+\frac{1}{4}x_3^2 \\ 
&&\textstyle    x_1x_8-x_2x_7+ \frac{1}{2}x_3x_5\\  
&&\textstyle  x_1x_{11}-x_2x_9+\frac{1}{4}x_5^2 \\
&&\textstyle  x_1x_{10}-\frac{1}{2}x_3x_7+x_4x_5 \\
&&\textstyle x_1x_{12}-x_3x_9+\frac{1}{2}x_5x_7 \\
&&\textstyle x_1x_{13}-x_4x_9+\frac{1}{4}x_7^2\\ 
&&\textstyle x_1x_{14}-\frac{1}{4}x_5x_{10}+\frac{1}{4}x_7x_8-x_6x_9 \\
&&\textstyle x_2x_{10}-\frac{1}{2}x_3x_8+x_5x_6\end{eqnarray} \begin{eqnarray}  
&&\textstyle x_2x_{12}-x_3x_{11}+\frac{1}{2}x_5x_8\\ 
&&\textstyle x_2x_{13}-\frac{1}{4}x_3x_{12}+\frac{1}{4}x_7x_8-x_6x_9 \\
&&\textstyle x_2x_{14}-x_6x_{11}+\frac{1}{4}x_8^2\\ 
&&\textstyle x_3x_{10}-2x_4x_8+2x_6x_7\\ 
&&\textstyle x_3x_{12}+x_5x_{10}-4x_4x_{11}+4x_6x_9\\ 
&&\textstyle x_3x_{13}-x_4x_{12}+\frac{1}{2}x_7x_{10}\\ 
&&\textstyle x_3x_{14}-x_6x_{12}+\frac{1}{2}x_8x_{10}\\ 
&&\textstyle x_5x_{12}-2x_7x_{11}+2x_8x_9\\ 
&&\textstyle x_5x_{13}-\frac{1}{2}x_7x_{12}+x_9x_{10}\\ 
&&\textstyle x_5x_{14}-\frac{1}{2}x_8x_{12}+x_{10}x_{11}\\ 
&&\textstyle x_4x_{14}-x_6x_{13}+\frac{1}{4}x_{10}^2\\ 
&&\textstyle x_7x_{14}-x_8 x_{13}+\frac{1}{2}x_{10}x_{12}\\  
&&\textstyle x_9x_{14}-x_{11}x_{13}+\frac{1}{4}x_{12}^2
\end{eqnarray}
\section{Equations of partial $A_3$ flag variety ${\rm FL}_{1,3}$ }\label{eqa3cd9}  Following the description given in \cite{rudakov}, we have the following decomposition of the second symmetric power of the contragradient representation \(V_{\lam}^*\). \[S^2V_{\lam}^*=V_1\oplus V_{e_{1}} \oplus V_{e_2} \oplus V_{e_3}.\] The vector spaces \(V_1, V_{e1},V_{e_2}\text{ and }V_{e_3}\) are 84, 20, 15 and 1 dimensional respectively. The defining equations of \(\fl_{1,3}\), are the basis of the linear subspaces of dimension 20, 15 and 1 of \(S^2V_{\lam}^*\). \[I=\left<Q\right>=\left< V_{e_1}\right> \cup\left < V_{e_2}\right> \cup \left<V_{e_3} \right>\subset S^2V_{\lam}^*.\] We compute these quadratic equations by using the GAP4 code  of appendix  of \cite{qs}  and they are listed below. 
 
\begin{eqnarray} \label{eqfl13}
&& \textstyle x_{1}x_{6}-x_{2}x_{3} \\
 && \textstyle    x_{1}x_{9}-x_{1}x_{7}-x_{1}x_{8}-x_{3}x_{4}+x_{2}x_{5} \\ 
&& \textstyle  x_{1}x_{12}-x_{3}x_{7}-x_{3}x_{8}+
x_{5}x_{6} \\ && \textstyle  x_{1}x_{11}+x_{2}x_{9}-x_{2}x_{8}-x_{4}x_{6} \\ && \textstyle  x_{1}x_{10}-x_{5}x_{4} \\ && \textstyle  
  x_{1}x_{14}+x_{3}x_{10}-x_{5}x_{9} \\ && \textstyle 
  x_{1}x_{13}-x_{2}x_{10}+x_{4}x_{7} \\ && \textstyle  
  x_{1}x_{15}+x_{5}x_{11}-x_{4}x_{12}-x_{6}x_{10}+x_{9}x_{7}+x_{9}x_{8}-x_{7}x_{8}-
x_{8}^2 \\ && \textstyle  x_{3}x_{11}+x_{2}x_{12}-x_{6}x_{8} \\ && \textstyle  x_{3}x_{14}-x_{5}x_{12}\end{eqnarray} \begin{eqnarray} \\ && \textstyle  x_{3}x_{13}-x_{2}x_{14}
    -x_{5}x_{11}+x_{4}x_{12}-x_{9}x_{8}+x_{7}x_{8}+x_{8}^2 \\ && \textstyle  x_{3}x_{15}-x_{6}x_{14}+x_{7}x_{12} \\ && \textstyle 
 x_{2}x_{13}-x_{4}x_{11} \\ && \textstyle  x_{2}x_{15}+x_{6}x_{13}-x_{9}x_{11} \\ && \textstyle  x_{5}x_{13}+x_{4}x_{14}-x_{8}x_{10} \\ && \textstyle  
  x_{5}x_{15}+x_{9}x_{14}-x_{8}x_{14}-x_{12}x_{10} \\
   && \textstyle  x_{4}x_{15}-x_{7}x_{13}-x_{8}x_{13}+
    x_{11}x_{10} \\ && \textstyle  x_{6}x_{15}-x_{12}x_{11} \\ && \textstyle  
  x_{9}x_{15}-x_{7}x_{15}-x_{8}x_{15}-x_{12}x_{13}+x_{11}x_{14} \\ && \textstyle  x_{10}x_{15}-x_{14}x_{13}  \\ && \textstyle  
 x_{1}x_{9}+x_{1}x_{7}-2x_{3}x_{4}-2x_{2}x_{5} \\ && \textstyle  
  x_{1}x_{12}-\frac{1}{2}x_{3}x_{9}+\frac{1}{2}x_{3}x_{7}-x_{5}x_{6} \\ && \textstyle  x_{1}x_{11}-\frac{1}{2}x_{2}x_{9}+\frac{1}{2}x_{2}x_{7}+
x_{4}x_{6} \\ && \textstyle  x_{1}x_{14}-x_{3}x_{10}+\frac{1}{2}x_{5}x_{9}
    -\frac{1}{2}x_{5}x_{7}-x_{5}x_{8} \\ && \textstyle x_{1}x_{13}+x_{2}x_{10}+\frac{1}{2}x_{4}x_{9}-\frac{1}{2}x_{4}x_{7}-x_{4}x_{8} \\ && \textstyle  x_{1}x_{15}-x_{3}x_{13}-
x_{2}x_{14}+x_{6}x_{10}-\frac{1}{2}x_{9}x_{7}+\frac{1}{2}x_{7}^2+x_{7}x_{8} \\ && \textstyle  x_{3}x_{11}-x_{2}x_{12}+
\frac{1}{2}x_{6}x_{9}+\frac{1}{2}x_{6}x_{7} \\ && \textstyle  x_{3}x_{13}+x_{2}x_{14}-x_{5}x_{11}-x_{4}x_{12}+\frac{1}{2}x_{9}^2
-\frac{1}{2}x_{9}x_{8}-\frac{1}{2}x_{7}^2-\frac{1}{2}x_{7}x_{8} \\ && \textstyle  
  x_{3}x_{15}+x_{6}x_{14}+\frac{1}{2}x_{9}x_{12}-\frac{1}{2}x_{7}x_{12}-x_{8}x_{12} \\ && \textstyle  x_{2}x_{15}-
x_{6}x_{13}+\frac{1}{2}x_{9}x_{11}-\frac{1}{2}x_{7}x_{11}-x_{8}x_{11} \\ && \textstyle  
  x_{5}x_{11}+x_{4}x_{12}-\frac{1}{4}x_{9}^2+\frac{1}{4}x_{7}^2 \\&& \textstyle  x_{5}x_{13}-x_{4}x_{14}+
\frac{1}{2}x_{9}x_{10}+\frac{1}{2}x_{7}x_{10} \\ && \textstyle  x_{5}x_{15}-\frac{1}{2}x_{9}x_{14}+\frac{1}{2}x_{7}x_{14}+x_{12}x_{10} \\ && \textstyle  
  x_{4}x_{15}-\frac{1}{2}x_{9}x_{13}+\frac{1}{2}x_{7}x_{13}-x_{11}x_{10} \\ && \textstyle  x_{9}x_{15}+x_{7}x_{15}-2x_{12}x_{13}-
2x_{11}x_{14} \\ && \textstyle x_{1}x_{15}-x_{3}x_{13}+x_{2}x_{14}-x_{5}x_{11}+
    x_{4}x_{12}-x_{6}x_{10}-\frac{3}{8}(x_{9}^2+x_{7}^2)\notag\\&& +\frac{1}{4}x_{9}x_{7}+        \frac{1}{2}(x_{9}x_{8}-
x_{7}x_{8}-x_{8}^2)
\end{eqnarray}

\bibliographystyle{amsplain}
\bibliography{imran}

{\sc Lums School of Science and Engineering}\\{\sc U-Block, DHA, Lahore, Pakistan.}
\\[2mm]
{\it Email address}: {\tt i.qureshi@maths.oxon.org}

\end{document}